\documentclass[reqno]{amsart}
\usepackage{amssymb}
\usepackage{graphicx}
\usepackage[left]{lineno}
\usepackage{enumitem}
\usepackage{xcolor}

\usepackage{csquotes}

\newtheorem{theorem}{Theorem}[section]
\newtheorem*{theorem*}{Theorem}

\newtheorem{lemma}[theorem]{Lemma}

\newtheorem{proposition}[theorem]{Proposition}
\newtheorem{corollary}[theorem]{Corollary}

\theoremstyle{definition}
\newtheorem{definition}[theorem]{Definition}

\theoremstyle{remark}
\newtheorem{remark}[theorem]{Remark}
\numberwithin{equation}{section}

\usepackage[margin=2.54cm]{geometry}
\DeclareMathOperator*{\argmin}{arg\,min}
\DeclareMathOperator*{\dom}{dom}

\title{Metric gradient flows and variational principles}

	\author[]{Alberto Dom\'inguez Corella}
    
	\address{Institut f\"{u}r Stochastik und Wirtschaftsmathematik, VADOR E105-04, TU Wien, Wiedner Hauptstra{\ss}e 8, A-1040 Wien, \"Osterreich}

    \address{Institut f\"{u}r Mathematik und Wissenschaftliches Rechnen, Universit\"{a}t Graz, Heinrichstra{\ss}e 36, A-8010 Graz, \"Osterreich\\}

	\email{alberto.of.sonora@gmail.com}


\thanks{The author was supported by the Austrian Science Foundation (FWF) under grants P 36344-N and F 100800.}

\begin{document}

\begin{abstract}
We show that Ekeland-type variational principles can be derived from steepest descent curves in metric spaces. From this, we relate the existence of the Calabi flow to the validity of the Ekeland variational principle for the $K$-energy on the space of K\"ahler potentials.
\end{abstract}

\maketitle
\setcounter{tocdepth}{1}
\tableofcontents
\section{Introduction and background}
\subsection{Background, motivation and overview}
The now-celebrated Ekeland variational principle can be regarded as the core of modern variational analysis, as it has found quite diverse applications in the field. In a few words, it states that if a point nearly minimizes a functional, then it is possible to find a nearby point that still nearly minimizes it and is also \textit{almost critical}.
\begin{theorem*}[Ekeland {\cite[Theorem 1.1]{E_1974}}]
    Let $(X,d)$ be a complete metric space and $\Phi:X\to\mathbb R\cup\{+\infty\}$ a proper lower semicontinuous functional.  
	Let $u\in X$ and $\varepsilon>0$ be such that 
	\begin{align*}
		\Phi(u) \le  \inf\Phi + \varepsilon.
	\end{align*}
	For any $\lambda>0$ there exists $\hat u\in X$ such that 
	\begin{align*}
		(a)\,\, d(u,\hat u)\le\lambda\qquad\,\,\,\, (b)\,\,\Phi(\hat{u}) + \frac{\varepsilon}{\lambda} d(u,\hat u) \leq \Phi(u)\qquad \,\,\,\,(c)\,\,   \Phi(\hat u) \le \Phi(v) + \frac{\varepsilon}{\lambda} d(v,\hat u)\quad\forall v\in X.
	\end{align*}  
\end{theorem*}
\smallbreak
\noindent
The  variational principle has been extended in many directions, but rarely in ways that dispense with completeness (at the expense of other assumptions) or provide constructive arguments (at least in some particular cases). Most, if not all, extension proofs rely on the same set-theoretic order arguments as the original result. We discuss the literature at the end of the section.
\smallbreak\noindent
This manuscript is devoted to the derivation of variational principles based solely on the global existence of metric gradient flows. Although the conclusions of the Ekeland variational principle — proximity, energy decrease, and small slope — bear a close heuristic resemblance to their qualitative behavior, no existing work has established variational  principles directly from gradient-flow dynamics.
\smallbreak\noindent
There are two main motivations for this derivation. The first is that there are several known results whose proofs can be carried out using essentially the same techniques — either using gradient flows or the Ekeland variational principle — as, for example, the following.
\begin{itemize}[label=\raisebox{0.2ex}{\scriptsize\(\bullet\)}]
    \item (\textit{Existence of minimizers via Palais--Smale conditions}). The use of Palais--Smale conditions to guarantee the existence of minimizers originated in \cite[Theorem~1]{PS_1964}, initially via gradient-flow arguments, and later became a standard application of the Ekeland variational principle; see \cite[pp.~340--341]{E_1974}.
    \smallbreak\noindent
     \item (\textit{Error bounds and \L ojasiewicz inequalities}). The \L ojasiewicz inequality is known to imply a local error bound estimate; such an implication was obtained in \cite[Theorem~5]{BNPS_2017} (using gradient-flow arguments) and also in \cite[Theorem~2]{VT_2009} (using the Ekeland variational principle).
     \smallbreak\noindent
     \item (\textit{Minimax critical value theorems}). Results yielding critical values from minimax problems have been proved both by gradient-flow methods \cite[Theorem~12.5]{CH_1982} and by the Ekeland variational principle \cite[Theorem~(1.ter) and Corollary~(5)]{G_1991}.
     \smallbreak\noindent
     \item (\textit{Determination of convex functions via slopes}). Results determining convex functions up to an additive constant have been proved both via gradient-flow methods \cite[Theorem~3.1 and Corollary~3.1]{PSV_2021} and by employing the Ekeland variational principle \cite[Theorem~4.3]{IZ_2025}.
\end{itemize}
\noindent
These examples indicate that gradient flow techniques and variational principles are somehow related. In practice, gradient flows provide more information than the variational principle. We show that this is because there are points on the flow curves that recover many of the properties possessed by the principle. Intuitively, this also suggests the heuristic conclusion that, in some situations (not all), one may use the principle as a substitute in proofs based on gradient-flow arguments.
\smallbreak\noindent
The second motivation is that such a derivation is particularly relevant in settings where the ambient space is not complete, but gradient flows are nevertheless available. In those situations, the classical variational principle cannot be applied in its usual form, whereas flows may still provide a canonical trade-off between proximity and descent. Examples of this are the following.
\begin{itemize}[label=\raisebox{0.2ex}{\scriptsize\(\bullet\)}]
  \item (\textit{Geometric flows}). There are relevant examples of flows coming from differential geometry where the underlying configuration space is non-complete. For instance, the mean curvature flow can be viewed as the $L^{2}$-gradient flow of the area functional on a space of smooth immersions \cite{E_2004,M_2011}, and the Calabi flow as the gradient flow of Mabuchi's $K$-energy on the space of K\"ahler potentials (with respect to the Mabuchi--Semmes--Donaldson metric) \cite{D_1999,M_1987}.
  \smallbreak\noindent
  \item (\textit{Wasserstein flows}). Several PDEs, such as the porous medium  and Fokker--Planck equations with no-flux boundary conditions, can be formulated as Wasserstein gradient flows on the space of probability densities over a bounded domain \cite{AG_2013,FG_2021}. It is a standard fact that such a space is not complete.
\end{itemize}\noindent
With these motivations in mind, we now describe the main contents of the manuscript.
\smallbreak\noindent
Our main result states that it is possible to construct points satisfying a variational principle along trajectories of metric gradient flows (Theorems \ref{thm_gvp} and \ref{thm_gvp_alt}). The gradient-flow framework we use is the now-classic one of \cite[Part I]{AGS_2008}, based on strong upper gradients. When the upper gradient is the metric slope, we obtain a local version of the Ekeland variational principle (Corollary \ref{cor_lms}); we also give an instance in which the full principle can be recovered (Corollary \ref{cor_gms}).  
As an application, we relate the existence of the Calabi flow to the validity of the Ekeland variational principle for Mabuchi's $K$-energy on the space of K\"ahler potentials endowed with the Mabuchi--Semmes--Donaldson metric (despite its lack of completeness). From the derived principles, several standard consequences are obtained, including a fixed point theorem (Proposition \ref{P3}) and a version of the  metric deformation lemma (Proposition \ref{P2}).
\smallbreak\noindent
An advantage of more explicit constructions is that they allow for finer arguments when applying the variational principle. They may also provide a useful starting point for studying constructive or computational variants of such principles, although we do not pursue this question here.

\subsection{On the related literature} 
It has been more than half a century since the  theorem of Ekeland \cite[Theorem 1.1]{E_1974}. Its proof came from the order-theoretic ideas behind the seminal Bishop--Phelps lemma \cite[Lemma 1 and Theorem 2]{BP_1963} and the Br{\o}ndsted--Rockafellar theorem \cite[Lemma on p.~608]{BR_1965}. Since then, numerous extensions and refinements have appeared.  
Quoting \cite[Introduction]{KPS_2016},
\begin{quote}
The number of publications containing ``Ekeland variational principle''
in their title has exceeded 200.
\end{quote}
This can be easily verified through a standard MathSciNet/ZbMATH search. Since reviewing all the references extending or modifying the variational principle is beyond the scope of this paper, we provide a brief list of those consulted in the preparation of this manuscript.
\begin{itemize}[label=\raisebox{0.2ex}{\scriptsize\(\bullet\)}]
     \item (\textit{General perturbations}). Among the most influential refinements of the  principle are the Stegall principle \cite[Theorem on p.~172]{S_1978} (linear perturbations), the Borwein--Preiss principle \cite[Theorem 2.6]{BP_1987} (quadratic-type perturbations), and the smooth variational principle in \cite[Theorem 2.3]{DGZ_1993}; see also \cite{BZ_2010,DP_2009} for related principles.
     \smallbreak\noindent
     \item (\textit{Ambient spaces}). Versions of the principle beyond the classical metric setting include asymmetric locally convex spaces \cite{C_2012}; quasi-metric and quasi-pseudometric frameworks \cite{C_2025,C_2011,C_2019,C_2023}; and partial-metric and quasi-partial-metric settings \cite{AKV_2015,MTK_2017}. Extensions to locally complete (or sequentially lower complete) spaces are developed in \cite{BGG_2007,BGGG_2015,BL_2014,HLL_2006,HQ_2015,JQ_2003}, and are often accompanied by fixed-point (e.g.\ Caristi-type) or best-proximity results \cite{C_2011,C_2019,IZZ_2015,MTK_2017,Z_2019}. 
     \smallbreak\noindent
    \item (\textit{The convex case}). The Brøndsted--Rockafellar theorem can be regarded as the Ekeland variational principle for convex functions. The theorem has been extended in many directions, including operator-theoretic \cite{B_2013c,IS_2012,MS_2008,S_2008b}
 and convex-analytic \cite{AO_2024,CR_2003,CHP_2018,FI_2024,L_2015,Z_2023}. In \cite[Theorem~4.1]{C_2023b}, an explicit construction based on the resolvent of the convex subdifferential was given; a similar construction appeared in the proof of \cite[Proposition~2.3]{L_2015}. This construction was recently extended, in \cite[Section~1.2]{DV_2026}, to geodesic spaces by using proximal points in lieu of the resolvent.

     \smallbreak\noindent
    \item (\textit{Relations to completeness}). A recurring theme is that completeness properties can be characterized \cite{S_1981, W_1977}, or at least detected via suitably weakened Ekeland-type principles on restricted function classes, notably lower bounded Lipschitz functions \cite{B_2012b}, and in related completeness classes such as UC  spaces \cite{B_2013}; see also \cite{K_2018} and the monograph \cite{M_2009}.
\end{itemize}\noindent
None of the aforementioned references relate the principle to gradient flows, and we are not aware of any other sources that do so. To the best of our knowledge, this paper is the first to derive Ekeland-type variational principles from gradient flows.
\smallbreak\noindent
Concerning the theory of metric gradient flows, we give the following account. 
\begin{itemize}[label=\raisebox{0.2ex}{\scriptsize\(\bullet\)}]
   \item (\textit{Minimizing movements and curves of maximal slope}). The theory of metric gradient flows originates in the minimizing-movement approach. While its proximal precursor is the implicit scheme in \cite{M_1965}, its first metric formulation appears in \cite{DMT_1980}, where evolutions are described directly in metric spaces through curves of maximal decrease; see also \cite{A_1995,MST_1989}. Stability aspects of the theory have been studied in \cite{CG_2012,MRS_2012,S_2011}. For a modern account of the theory we refer to \cite{AGS_2004,AGS_2008,HM_2019}.
    \smallbreak\noindent
        \item (\textit{Wasserstein gradient flows}).
    The central examples in the metric theory come from Wasserstein gradient flows. 
    This direction began with the introduction of convexity along transport geodesics \cite{M_1997} and the observation that the Fokker--Planck equation can be obtained as a minimizing movement for the free energy in the quadratic Wasserstein distance \cite{JKO_1998}. The geometrical interpretation of diffusion-type equations as gradient flows in Wasserstein spaces was further developed in \cite{O_2001}.  The theory was subsequently placed in the metric framework in \cite[Part II]{AGS_2008}. 
    Modern and accessible references include \cite{ABS_2024,FG_2021,S_2015,S_2017}.
    \smallbreak\noindent
    \item (\textit{Gradient flows in Hadamard spaces}). Another source of examples comes from nonpositively curved spaces. In \cite{M_1998}, gradient flows for convex functionals on nonpositively curved metric spaces were constructed via minimizing movements. This theory was further developed through proximal point methods, product
formulae, and nonlinear semigroups in Hadamard spaces
\cite{B_2013b,B_2014,B_2015,BR_2014,S_2012}.
    \smallbreak\noindent
    \item (\textit{Further extensions of the metric setting}). Further extensions consider gradient-flow structures induced by
transport-type distances adapted to the underlying dynamics.
A discrete Benamou--Brenier-type construction for finite irreducible
reversible Markov chains was shown
in \cite{M_2011b} to yield an entropy gradient-flow interpretation of
the Markov semigroup. Related discrete and nonlocal Markovian formulations were developed next in
\cite{CHLZ_2012,E_2014,M_2013}.

\end{itemize}

\section{Gradient flows and variational principles} \label{sec:mgf}
\noindent
We present in this section the results discussed in the introduction, together with some examples illustrating the sharpness of the assumptions. All proofs are postponed until the final subsection. Throughout the section, $(X, d)$ will denote a metric space and $\Phi:X\to\mathbb R\cup\{+\infty\}$ a proper functional bounded from below. 

\subsection{Metric gradient flows}\label{mgfs}
Let us begin by recalling the definition of metric  gradient flow; preliminaries on metric derivatives and (strong) upper gradients are given in the appendix. 

\begin{definition}\label{dgf}
	Let $g:X\to\mathbb R\cup\{+\infty\}$ be an upper gradient of $\Phi$. 
	We say that a locally absolutely continuous curve  $v:(0,+\infty)\to X$ is  a \textit{gradient flow} of  $\Phi$ with respect to $g$ if 
	\begin{enumerate}
		\item $v(t)\in\dom\Phi$ for all $t\in(0,+\infty)$; 
		\smallbreak 
		\item $\Phi(v(t)) - \Phi(v(s))  \le -\displaystyle\frac{1}{2} \displaystyle\int_{s}^t |\dot{v}|(r)^2\,dr  - \displaystyle\frac{1}{2} \displaystyle\int_{s}^t g(v(r))^2\,dr$ for all $s,t\in(0,+\infty)$.
	\end{enumerate}
\end{definition}
\noindent
We mention that Definition \ref{dgf} is equivalent to the usual one \cite[Definition 1.3.2]{AGS_2008}, see \cite[Remark 1.3.3]{AGS_2008}. 
We now state a simple but very useful result concerning the initial points of gradient flows. 
\begin{lemma}\label{AB}
    Let $g:X\to\mathbb R\cup\{+\infty\}$ be an upper gradient of $\Phi$. 
	Let   $v:(0,+\infty)\to X$ be a {gradient flow} of $\Phi$ with respect to $g$. Suppose there exists $u\in\dom\Phi$ such that  
    \[
    \lim_{t\longrightarrow 0^+} v(t) = u\quad \text{and}\quad   \lim_{t\longrightarrow 0^+}\Phi(v(t)) = \Phi(u). 
\]
Define $v(0):=u$. Then $v:[0,+\infty)\to X$ is locally absolutely continuous and 
\begin{align*}
    \Phi(v(t)) - \Phi(v(s))  = - \displaystyle\int_{s}^t |\dot{v}|(r)^2\,dr  =- \int_{s}^t g(v(r))^2\,dr
\end{align*}
 for all $s,t\in[0,+\infty)$. In particular, $|\dot v|\in L^2_{\operatorname{loc}}([0,+\infty))$ and $g\circ v\in L^2_{\operatorname{loc}}([0,+\infty))$.
\end{lemma}
\noindent
We recall that \(L^2_{\mathrm{loc}}([0,+\infty))\) is the space of measurable functions from
\([0,+\infty)\) into \(\mathbb{R}\) whose restriction to every compact subset
\(K \subset [0,+\infty)\) is square integrable. We note that Lemma \ref{AB} may fail
if one assumes only the existence of the limit of the trajectory as
\(t \longrightarrow 0^+\) without requiring convergence of the energy along the
trajectory; see Subsubsection \ref{refloj} below. 

\subsection{The variational principle}
In this subsection we formulate the main results of the paper. We recall that all proofs in this section are
deferred to the end of the section.
\smallbreak\noindent
Given an upper gradient $g$ of $\Phi$, we define $\mathcal M(g)$ as the set of all
$u\in \dom\Phi$ such that there exists a gradient flow
$v:(0,+\infty)\to X$ of $\Phi$ with respect to $g$ satisfying
\[
    \lim_{t\longrightarrow 0^+} v(t) = u
    \quad\text{and}\quad
    \lim_{t\longrightarrow 0^+}\Phi(v(t)) = \Phi(u).
\]
The set \(\mathcal M(g)\) consists of those points in the effective domain of
\(\Phi\) from which a gradient flow with well--defined initial position and
energy can be started.

\subsubsection{The general variational principle}
The first result is a   variational principle based solely on the existence of a metric gradient
flow with respect to a given upper gradient.
\begin{theorem}\label{thm_gvp}
Assume that there exists an upper gradient  $g:X\to\mathbb R\cup\{+\infty\}$ of $\Phi$ such that $\mathcal M(g)\neq \emptyset$. 
	Let $u\in\mathcal M(g)$ and $\varepsilon>0$ be such that 
	\begin{align*}
		\Phi(u) \le \inf\Phi + \varepsilon.
	\end{align*}
	For any $\lambda>0$ there exists $\hat u\in \mathcal M(g)$ such that 
	\begin{align*}
		(a)\,\,\, d(u,\hat u)\le\lambda\qquad\,\,\,\, (b)\,\,\,\Phi(\hat{u})  \leq \Phi(u)\qquad \,\,\,\,(c)\,\,\,   g(\hat u) \le \frac{\varepsilon}{\lambda}. 
	\end{align*} 
\end{theorem}\medskip\noindent
We now make some remarks concerning Theorem \ref{thm_gvp}.
\begin{remark}
The natural candidate for an upper gradient is usually the so-called
\emph{relaxed slope} \cite[Section~2.3]{AGS_2008}; under suitable hypotheses,
gradient flows with respect to it can be constructed by means of the
\emph{generalized minimizing movement} scheme, see \cite[Theorem~2.3.3]{AGS_2008}.
\end{remark}

\begin{remark}
    The theorem does not require completeness of the underlying metric space; the only substantial
hypothesis is the global existence of a gradient flow, which can be relaxed to the global
existence of $p$-gradient flows \cite[Definition 1.3.2]{AGS_2008}  with $p\in(1,+\infty)$ (the arguments in the proof
repeat under trivial adaptations).
\end{remark}

\begin{remark}
The element $\hat u\in\mathcal M(g)$ in the theorem is selected along a gradient
flow $v:(0,+\infty)\to X$ starting from $u\in\mathcal M(g)$; either $\hat u=u$, or
there exists $\tau>0$ such that $\hat u=v(\tau)$.
\end{remark}

\begin{remark}\label{strictimpro}
It follows from the construction in the proof of the theorem that the element
$\hat u\in\mathcal M(g)$ actually satisfies the strict estimate
$g(\hat u)<\varepsilon/\lambda$.
\end{remark}

\subsubsection{The variational principle under more regularity}
Under the additional assumption of lower semicontinuity of the upper gradient at the
reference point, the conclusion of Theorem~\ref{thm_gvp} can be strengthened. We present now the second result. 

  \begin{theorem}\label{thm_gvp_alt}
Assume that there exists an upper gradient  $g:X\to\mathbb R\cup\{+\infty\}$ of $\Phi$ such that  $\mathcal M(g)\neq \emptyset$.
	Let $u\in\mathcal M(g)$ and $\varepsilon>0$ be such that 
	\begin{align*}
		\Phi(u) \le \inf\Phi + \varepsilon\qquad\text{and}\qquad g(u)\le\liminf_{v\longrightarrow u} g(v).
	\end{align*}
   Then 
	for any $\lambda>0$ there exists $\hat u\in \mathcal M(g)$ such that 
	\begin{align*}
		(a)\,\,\, d(u,\hat u)\le\lambda\qquad\,\,\,\, (b)\,\,\,\Phi(\hat{u}) +\frac{\varepsilon}{\lambda} d(u,\hat u) \leq \Phi(u)\qquad \,\,\,\,(c)\,\,\,   g(\hat u) \le \frac{\varepsilon}{\lambda}. 
	\end{align*} 
\end{theorem}\medskip\noindent
We mention that the lower semicontinuity assumption on the upper gradient is not superfluous; the sole assumptions in Theorem \ref{thm_gvp} are not sufficient for the bound in Theorem \ref{thm_gvp_alt}-$(b)$, see Subsubsection \ref{exeke}. We also mention that in contrast with Theorem~\ref{thm_gvp} (see Remark \ref{strictimpro}), the estimate \(g(\hat u)\le \varepsilon/\lambda\) cannot be improved to a strict inequality in general; see Subsubsection~\ref{ex:sharp_c}.

\subsection{Ekeland-type variational principles}
In this subsection we consider specific upper gradients and specialize the variational principles of the previous subsection accordingly.
\subsubsection{A local version of the principle} 
The \textit{local metric slope} of $\Phi$ at a point $u\in X$ is given by 
\[
    |\nabla\Phi|(u):=
    \left\{
    \begin{array}{lcc}
    \displaystyle\max\left\{\limsup_{v\longrightarrow u}\frac{\Phi(u)-\Phi(v)}{d(u,v)},\hspace{0.02cm}0\right\} & \text{if} & u\in\dom\Phi \\ \\
    +\infty & \text{if} & u\notin\dom\Phi.
    \end{array}
    \right.
\]
This quantity is always nonnegative, and vanishes at local minimizers. If $X$ is a real normed space, and $\Phi$ a Fr\'echet differentiable functional, then $|\nabla\Phi|(u)=\|D\Phi(u)\|_{X^*}$ for all $u\in X$. This identity is one of the reasons for the term metric slope.
\smallbreak\noindent
From the definition of the metric slope and Theorem \ref{thm_gvp}, we obtain the following corollary.
\begin{corollary}\label{cor_lms}
Assume that the local metric slope $|\nabla\Phi|:X\to\mathbb R\cup\{+\infty\}$ is an  upper gradient of $\Phi$ and that 
$\mathcal M(|\nabla\Phi|)\neq\emptyset$. 
	Let $u\in\mathcal M(|\nabla\Phi|)$ and $\varepsilon>0$ be such that 
	\begin{align*}
		\Phi(u) \le  \inf\Phi + \varepsilon.
	\end{align*}
	For any $\lambda>0$ there exist $\hat u\in \mathcal M(|\nabla\Phi|)$ and a neighborhood $\mathcal V$ of $\hat u$ such that 
	\begin{align*}
		(a)\,\,\, d(u,\hat u)\le\lambda\qquad\,\,\,\, (b)\,\,\,\Phi(\hat{u})  \leq \Phi(u)\qquad \,\,\,\,(c)\,\,\,   \Phi(\hat u) < \Phi(v) + \frac{\varepsilon}{\lambda} d(v,\hat u)\quad\forall v\in \mathcal V\setminus\{\hat u\}. 
	\end{align*} 
\end{corollary}\medskip\noindent
We mention that, under more regularity of the metric slope (lower semicontinuity), one can improve the estimate in Corollary \ref{cor_lms}-$(b)$ in accordance with Theorem \ref{thm_gvp_alt}.

\subsubsection{A global version of the principle} 
The \textit{global metric slope} of $\Phi$ at a point $u\in X$ is given by 
\[
    |\partial\Phi|(u):=
    \left\{
    \begin{array}{lcc}
    \displaystyle\max\left\{\sup_{v\in X\setminus\{u\}}\frac{\Phi(u)-\Phi(v)}{d(u,v)},\hspace{0.02cm}0\right\} & \text{if} & u\in\dom\Phi \\ \\
    +\infty & \text{if} & u\notin\dom\Phi.
    \end{array}
    \right.
\]
This quantity is always nonnegative, and vanishes exactly at the global minimizers of $\Phi$. If $X$ is a real normed space and $\Phi$ a proper convex functional, then the metric slope coincides with the distance of the convex subdifferential to zero, i.e., $|\partial \Phi|(u) =\operatorname{dist}\big(0,\partial \Phi(u)\big)$ for all $u\in X$, where $\partial\Phi:X\rightrightarrows X^*$ denotes the convex subdifferential of $\Phi$; see \cite[Proposition 2.1 $(vii)$]{K_2015}. 
\smallbreak\noindent
A consequence of Theorem~\ref{thm_gvp_alt} is that, when gradient flows with respect to the global slope exist, the Ekeland variational principle can be recovered in its usual form.
\begin{corollary}\label{cor_gms}
Suppose that the functional $\Phi:X\to\mathbb R\cup\{+\infty\}$ is lower semicontinuous and that $\mathcal M(|\partial\Phi|)\neq\emptyset$.  
	Let $u\in\mathcal M(|\partial\Phi|)$ and $\varepsilon>0$ be such that 
	\begin{align*}
		\Phi(u) \le  \inf\Phi + \varepsilon.
	\end{align*}
	For any $\lambda>0$ there exists $\hat u\in \mathcal M(|\partial\Phi|)$ such that 
	\begin{align*}
		(a)\,\,\, d(u,\hat u)\le\lambda\qquad\,\,\,\, (b)\,\,\,\Phi(\hat{u}) +\frac{\varepsilon}{\lambda} d(u,\hat u) \leq \Phi(u)\qquad \,\,\,\,(c)\,\,\,   \Phi(\hat u) \le \Phi(v) + \frac{\varepsilon}{\lambda} d(v,\hat u)\quad\forall v\in X.
	\end{align*}  
\end{corollary}\medskip\noindent
We now make a couple of remarks on the previous corollary.

\begin{remark}
The lower semicontinuity of $\Phi$ guarantees that the global metric slope
$|\partial\Phi|$ is a (strong) upper gradient; see \cite[Theorem~1.2.5]{AGS_2008}
(the result does not require completeness of $X$).
\end{remark}

\begin{remark}
In many situations, the local and global metric slopes coincide (this property is referred to as the \emph{slope cone property} in \cite[Definition~3.1]{CG_2012}), e.g., when the functional is geodesically convex; see \cite[Definition 2.4.3]{AGS_2008} and \cite[Theorem~2.4.9]{AGS_2008}.
\end{remark}

\subsection{Examples and counterexamples}

\subsubsection{Local square-integrability of the metric derivative}\label{refloj}
Given an upper gradient $g:X\to\mathbb R\cup\{+\infty\}$, the set $\mathcal{M}(g)$ consists of those initial states from which a gradient flow is well defined both in position and in energy. This is essential in the proofs of Theorems~\ref{thm_gvp} and~\ref{thm_gvp_alt}, since Lemma~\ref{AB} yields the required local square--integrability estimates only under energy convergence. The following elementary example shows that trajectory convergence alone does not prevent the metric derivative from failing to be locally square integrable near time zero, even in Euclidean spaces.
\smallbreak\noindent
Let $X=\mathbb R^d$ with the Euclidean distance. Let $\Phi:X\to\mathbb R\cup\{+\infty\}$ and $g:X\to \mathbb R\cup\{+\infty\}$ be given by 
\[
\Phi(u) := \left\{
\begin{array}{lcc}
 \dfrac{1}{|u|} & \text{if} & u \in \mathbb R^d\setminus\{0\} \\ \\
0 & \text{if} & u = 0 
\end{array}
\right.
\qquad \text{and} \qquad
g(u) := \left\{
\begin{array}{lcc}
\dfrac{1}{|u|^{2}} & \text{if} & u \in \mathbb R^d\setminus\{0\} \\ \\
0 & \text{if} & u = 0.
\end{array}
\right.
\]
Fix a unit vector $e\in\mathbb R^d$ and consider the curve $v:(0,+\infty)\to X$ given by 
\[
v(t):=\sqrt[3]{3t}\,e.
\]
It is easy to check that this curve is a gradient flow of $\Phi$ with respect to $g$ and that $\displaystyle\lim_{t\longrightarrow0^+}v(t)=0$, but  $|\dot v|\notin L^{2}_{\mathrm{loc}}([0,+\infty))$.

\subsubsection{Necessity of lower semicontinuity of the upper gradient in Theorem \ref{thm_gvp_alt}}\label{exeke}
The improvement from Theorem~\ref{thm_gvp} to Theorem~\ref{thm_gvp_alt} relies on an additional assumption at the reference point, namely a lower semicontinuity condition for the upper gradient there. The next example shows that this assumption cannot, in general, be omitted.
\smallbreak\noindent
Let $X=\bigl([0,1]\times\{0\}\bigr)\cup\bigl(\{0\}\times[0,1]\bigr)$ 
equipped with the intrinsic path metric, i.e., 
\[
d(u,v)=
\begin{cases}
|u_1-v_1| & \text{if } u=(u_1,0)\ \text{and}\ v=(v_1,0)\\[2mm]
|u_2-v_2| & \text{if } u=(0,u_2)\ \text{and}\ v=(0,v_2)\\[2mm]
u_1+v_2   & \text{if } u=(u_1,0)\ \text{and}\ v=(0,v_2)\\[2mm]
u_2+v_1   & \text{if } u=(0,u_2)\ \text{and}\ v=(v_1,0).
\end{cases}
\]
Consider \(\Phi:X\to\mathbb R\cup\{+\infty\}\) and \(g:X\to\mathbb R\cup\{+\infty\}\) given by
\[
\Phi(u)=
\begin{cases}
-u_1 & \text{if } u=(u_1,0) \ \text{with}\ u_1\in[0,1]\\[2mm]
-3u_2 & \text{if } u=(0,u_2) \ \text{with}\ u_2\in(0,1]
\end{cases}
\qquad\text{and}\qquad
g(u)=
\begin{cases}
0 & \text{if } u=(1,0)\\[2mm]
1 & \text{if } u=(u_1,0) \ \text{with}\ u_1\in(0,1)\\[2mm]
3 & \text{if } u=(0,u_2) \ \text{with}\ u_2\in[0,1].
\end{cases}
\]
We note that \(g\) is not lower semicontinuous at $(0,0)$. Consider the curve \(v:(0,+\infty)\to X\) given by $v(t)=(\min\{t,1\},0)$. It is not difficult to see that $v$ is a gradient flow of $\Phi$ with respect to $g$ and hence that $(0,0)$ belongs to $\mathcal M(g)$. 
\smallbreak\noindent
To see that Theorem~\ref{thm_gvp_alt}-$(b)$ can fail without lower semicontinuity, take $u=(0,0)$, \(\varepsilon=3\) and \(\lambda=3/2\). Then \(\Phi(u)=0=\inf_X\Phi+\varepsilon\) and \(\varepsilon/\lambda=2\). If \(\hat u\in X\) satisfies items $(a)$ and $(c)$ of the theorem, i.e.,
\[
d(u,\hat u)\le \lambda\quad\text{and}\quad g(\hat u)\le 2,
\]
then necessarily \(\hat u=(\hat u_1,0)\) for some \(\hat u_1\in(0,1]\) since \(g\equiv3\) on \(\{0\}\times[0,1]\). For such \(\hat u\), \(d(u,\hat u)=\hat u_1\) and \(\Phi(\hat u)=-\hat u_1\), so $(b)$ would read
\[
\Phi(\hat u)+\frac{\varepsilon}{\lambda}d(u,\hat u)\le \Phi(u)
\quad\Longleftrightarrow\quad
-\hat u_1+2\hat u_1\le 0,
\]
yielding a contradiction since $\hat u_1>0$.

\subsubsection{Sharpness of bound in Theorem \ref{thm_gvp_alt}-$(c)$}\label{ex:sharp_c}
Even under the additional regularity hypothesis required in Theorem~\ref{thm_gvp_alt}, one cannot, in general, strengthen the conclusion by requiring the selected point to have an upper gradient strictly below the prescribed threshold. The following example in Euclidean space shows that this can already happen for a smooth  strictly convex energy.
\smallbreak\noindent
Let $X=\mathbb{R}^d$ with the Euclidean distance. Define $\Phi:X\to\mathbb{R}\cup\{+\infty\}$ and $g:X\to\mathbb{R}\cup\{+\infty\}$ by
\[
\Phi(x)=\frac12 |x|^2
\qquad\text{and}\qquad
g(x)=|x|.
\]
 Fix a unit vector $u\in\mathbb{R}^d$. It is not difficult to verify that the curve $v:(0,+\infty)\to X$ given by $v(t):=u e^{-t}$ is a gradient flow of $\Phi$ with respect to $g$ and that $u\in\mathcal M(g)$. 
\smallbreak\noindent
We now show that the bound in item~(c) of Theorem~\ref{thm_gvp_alt} cannot, in general, be improved to a strict inequality. Since $\inf\Phi=0$ and $\Phi(u)=1/2$, choosing $\varepsilon=1/2$ gives 
\[
\Phi(u)\le \inf\Phi+\varepsilon.
\]
Choose $\lambda=1/2$, so that $\varepsilon/\lambda=1$ and $g(u)=1$.
Suppose that $\hat u\in X$ satisfies items $(a)$ and $(b)$ of Theorem~\ref{thm_gvp_alt}, i.e.,
\[
|\hat u-u|\le\lambda
\quad\text{and}\quad
\Phi(\hat u)+\frac{\varepsilon}{\lambda}\,d(u,\hat u)\le\Phi(u).
\]
This can be equivalently rewritten as 
\begin{align}
    \frac12 |\hat u|^2 + |\hat u-u| \le \frac12.
\label{soin}
\end{align}
 By the reverse triangle inequality, $|\hat u|\ge 1-|\hat u-u|$, hence
\[
\frac12|\hat u|^2+|\hat u-u| \ge \frac12(1-|\hat u-u|)^2+|\hat u-u|=\frac12+\frac12|\hat u-u|^2.
\]
Together with (\ref{soin}), this yields $|\hat u-u|=0$ and thus $\hat u=u$. For this unique admissible point,
\[
g(\hat u)=g(u)=1=\frac{\varepsilon}{\lambda},
\]
so one cannot, in general, replace (c) by the strict inequality $g(\hat u)<\varepsilon/\lambda$.

\subsection{Proofs}
\subsubsection{Proof of Lemma \ref{AB}}
   Let $T>0$ be arbitrary. From \cite[Remark 1.3.3]{AGS_2008}, 
    \[
          \Phi(v(T)) - \Phi(v(s))  = - \displaystyle\int_{s}^T |\dot{v}|(r)^2\,dr  =- \int_{s}^T g(v(r))^2\,dr\quad \forall s\in(0,T].
    \]
    Taking limit as $s\longrightarrow0^+$, and using that  $\lim_{s\longrightarrow 0^+}\Phi(v(s)) = \Phi(u)$, we obtain that both $|\dot v|$ and $g\circ v$ belong to $L^2([0,T])$. Since $T>0$ was arbitrary, the result follows. 

\subsubsection{Proof of Theorem \ref{thm_gvp}}
Let $u\in\mathcal M(g)$ and $\varepsilon>0$ satisfy $\Phi(u)\le\inf\Phi + \varepsilon$.  
Since $u\in\mathcal M(g)$, there exists   a gradient flow  
 $v:(0,+\infty)\to X$ of $\Phi$ with respect to $g$ satisfying
\begin{align*}
    \lim_{t\longrightarrow 0^+} v(t) = u\quad \text{and}\quad   \lim_{t\longrightarrow 0^+}\Phi(v(t)) = \Phi(u). 
\end{align*}  
 Due to Lemma \ref{AB}, $|\dot v|\in L^2_{\operatorname{loc}}([0,+\infty))$, $g\circ v\in L^2_{\operatorname{loc}}([0,+\infty))$ and 
\begin{align*}
    \Phi(v(t)) - \Phi(v(s))  = - \displaystyle\int_{s}^t |\dot{v}|(r)^2\,dr  =- \int_{s}^t g(v(r))^2\,dr
\end{align*}
 for all $s,t\in[0,+\infty)$. 
Let $\lambda>0$ be given and define
	\begin{align*}
		\tau:=\inf \big\{t>0:\,\, g(v(t))<\varepsilon/\lambda \big\}. 
	\end{align*}
We now divide the proof into two cases.
\smallbreak\smallbreak\noindent
(\textbf{Case $\tau=0$}). By definition of infimum, there exists a sequence $(t_n)_{n\in\mathbb N}$ converging to zero such that $g(v(t_n))< \varepsilon/\lambda$. Choose $N\in\mathbb N$ large enough so that $d(v(t_N),u)\le \lambda$. Then, $\hat u:=v(t_N)$ satisfies items $(a)$-$(c)$ and belongs to $\mathcal M(g)$.
\smallbreak\noindent
(\textbf{Case} $\tau >0$ and $\Phi(v(\tau))>\inf\Phi$).
 Observe that  $g(v(t))\ge\varepsilon/\lambda$ for all $t\in(0,\tau)$ and hence
	\begin{align*}
		\varepsilon \ge \Phi(u)-\inf\Phi> \Phi(u) - \Phi(v(\tau)) = \int_{0}^{\tau} g(v(t))^2\,dt\ge \tau\frac{\varepsilon^2}{\lambda^2}.
	\end{align*} 
	We conclude that $\lambda^2/\varepsilon> \tau$. 
    Choose $\theta\in(0,1)$ such that 
    \begin{align}\label{thetaroom}
        \frac{\theta^2\lambda^2}{\varepsilon}>\tau.
    \end{align}
    From (\ref{thetaroom}), we obtain   
        \begin{equation}\label{gen2}
            \begin{aligned}
            d(v(\tau),u)
            &\le \int_{0}^{\tau}|\dot v|(s)\, ds \le \sqrt{\tau} \sqrt{\int_{0}^{\tau} |\dot v|(s)^2\, ds}= \sqrt{\tau}\sqrt{\Phi(u)-\Phi(v(\tau))}< \sqrt{\tau}\sqrt{\varepsilon} < \theta\lambda.
            \end{aligned}
        \end{equation}
    Due to the absolute continuity of the integral, there exists $\delta>0$ such that 
\begin{align}\label{absconint}
    \int_{\tau}^{\tau+\delta} |\dot v|(s)\, ds\le (1-\theta)\lambda.
\end{align}
    By definition of infimum, there exists  $\tau_\delta\in[\tau,\tau+\delta)$ such that $g(v(\tau_\delta))<\varepsilon/\lambda$. We prove that $\hat u:=v(\tau_\delta)$ satisfies the desired properties.
    Combining (\ref{gen2}) with (\ref{absconint}) and using the triangle inequality,
    \[
        d(\hat u,u)\le d(v(\tau_\delta),v(\tau)) + d(v(\tau),u) < \int_{\tau}^{\tau+\delta} |\dot v|(t)\,dt + \theta\lambda \le \lambda. 
    \]
   This shows the validity of item $(a)$; the one of $(c)$ follows directly from the definition of $\tau_\delta$. Item $(b)$ holds due to the monotonicity of $\Phi\circ v$. By construction, clearly $\hat u=v(\tau_{\delta})$ belongs to $\mathcal M(g)$.
    \smallbreak\noindent
    (\textbf{Case} $\tau>0$ and $\Phi(v(\tau))=\inf\Phi$).
Since $\Phi\circ v$ is nonincreasing,
\[
  \Phi(v(t)) = \inf\Phi \qquad\forall t\in[\tau,+\infty).
\]
By Lemma~\ref{AB}, we have
\[
  0=\Phi(v(T))-\Phi(v(\tau)) 
  = -\int_{\tau}^T|\dot v|^2(t)\,dt 
  = -\int_{\tau}^T g(v(t))^2\,dt\quad\forall T\in[\tau, +\infty).
\]
It follows that $|\dot v|(t)=0$ and $g(v(t))=0$ for a.e.\ $t\in[\tau,+\infty)$.
Thus there exists $\hat u\in X$ with $g(\hat u)=0$ such that $v(t)=\hat u$ for all
$t\ge\tau$. 
By Cauchy--Schwarz,
\[
  d(u,\hat u)
  \le \int_0^\tau |\dot v|(t)\,dt
  \le \sqrt{\tau}\Bigl(\int_0^\tau|\dot v|^2(t)\,dt\Bigr)^{1/2}
  = \sqrt{\tau}\,\sqrt{\Phi(u)-\Phi(\hat u)}
  \le \sqrt{\tau}\,\sqrt{\varepsilon}.
\]
On the other hand, since $g(v(t))\ge\varepsilon/\lambda$ for all $t\in(0,\tau)$,
\[
  \Phi(u)-\Phi(\hat u)
  = \int_0^\tau g(v(t))^2\,dt
  \ge \tau\,\frac{\varepsilon^2}{\lambda^2}.
\]
From this $\tau\le\lambda^2/\varepsilon$ and hence $ d(u,\hat u)\le \sqrt{\tau}\,\sqrt{\varepsilon}\le \lambda$.
We see that $\hat u$ satisfies all required properties. \hfill\(\square\)

\subsubsection{Proof of Theorem \ref{thm_gvp_alt}}
Since $u\in\mathcal M(g)$, there exists   a gradient flow  
 $v:(0,+\infty)\to X$ of $\Phi$ with respect to $g$ satisfying
\begin{align*}
    \lim_{t\longrightarrow 0^+} v(t) = u\quad \text{and}\quad   \lim_{t\longrightarrow 0^+}\Phi(v(t)) = \Phi(u). 
\end{align*}  
Define $\tau:=\inf \big\{t>0:\,\, g(v(t))\le\varepsilon/\lambda \big\}.$
We now divide the proof into two cases.
\smallbreak\noindent
(\textbf{Case $\tau=0$}). There exists a sequence $(t_n)_{n\in\mathbb N}$ converging to zero such that $g(v(t_n))\le \varepsilon/\lambda$. Let $\hat u:=u$ and observe that 
\[
    g(\hat u)=g(u)\le \liminf_{v\longrightarrow u} g(v) \le \liminf_{n\longrightarrow +\infty} g(v(t_n))\le \varepsilon/\lambda.
\]
We see that $\hat u$ satisfies trivially all items.
\smallbreak\noindent
(\textbf{Case $\tau >0$}). 
 Observe that  $g(v(t))>\varepsilon/\lambda$ for all $t\in(0,\tau)$ and hence
	\begin{align}\label{gen1l}
		\varepsilon \ge \Phi(u) - \Phi(v(\tau)) = \int_{0}^{\tau} g(v(t))^2\,dt> \tau\frac{\varepsilon^2}{\lambda^2}.
	\end{align} 
	We conclude that  $\lambda^2/\varepsilon> \tau$. Choose $\theta\in(0,1)$ such that $\theta^2\lambda^2/\varepsilon> \tau$.
    From this, 
	\begin{align*}
		d(v(\tau),u)\le \int_{0}^{\tau}|\dot v|(s)\, ds \le \sqrt{\tau} \sqrt{\int_{0}^{\tau} |\dot v|(s)^2\, ds} < \frac{\theta\lambda}{\sqrt{\varepsilon}} \sqrt{\Phi(u)-\Phi(v(\tau))} \le\theta \lambda. 
	\end{align*}
    Due to the absolute continuity of the integral, there exists $\delta_1>0$ such that 
\begin{align}\label{absconintl}
    \int_{\tau}^{\tau+\delta_1} |\dot v|(s)\, ds\le (1-\theta)\lambda.
\end{align}
     From (\ref{gen1l}) and the monotonicity of $\Phi\circ v$, we see that we can find $\delta_2>0$ such that 
     \begin{align}\label{gen3l}
         \Phi(u) - \Phi(v(s)) \ge \Phi(u) - \Phi(v(\tau)) > (\tau+\delta_2)\frac{\varepsilon^2}{\lambda^2}\quad\forall s\in[\tau,+\infty).
     \end{align}
    Define $\delta:=\min\{\delta_1,\delta_2\}$. 
    By definition of infimum, there exists  $\tau_\delta\in[\tau,\tau+\delta)$ such that $g(v(\tau_\delta))\le \varepsilon/\lambda$. We will prove that $\hat u:=v(\tau_\delta)$ satisfies the desired properties.
    Using the triangle inequality and (\ref{absconintl}), 
    \[
        d(\hat u,u)\le d(v(\tau_\delta),v(\tau)) + d(v(\tau),u) < \int_{\tau}^{\tau+\delta} |\dot v|(s)\,ds + \theta\lambda \le \lambda. 
    \]
   This shows the validity of item $(a)$; the one of $(c)$ follows directly from the definition of $\tau_\delta$. For the proof of item $(b)$, observe that, using (\ref{gen3l}),
\begin{align*}
  \frac{\varepsilon}{\lambda} d(u,v(\tau_\delta)) \le\frac{\varepsilon}{\lambda}  \sqrt{\tau_\delta} \sqrt{\Phi(u)-\Phi(v(\tau_\delta))}< \sqrt{\Phi(u)-\Phi(v(\tau_\delta))}\sqrt{\Phi(u)-\Phi(v(\tau_\delta))}={\Phi(u)-\Phi(v(\tau_\delta))}.
\end{align*}
This completes the proof. \hfill\(\square\) 
\subsubsection{Proof of Corollary \ref{cor_lms}}
By Theorem \ref{thm_gvp} and Remark \ref{strictimpro}, one can find $\hat u\in\mathcal M(|\nabla\Phi|)$ such that $d(u,\hat u)\le\lambda$, $\Phi(\hat u)\le\Phi(u)$ and $|\nabla\Phi|(\hat u)<\varepsilon/\lambda.$
By definition of local metric slope, 
\[
    \frac\varepsilon\lambda>|\nabla\Phi|(\hat u)\ge\limsup_{v\to \hat u}  \frac{\Phi(\hat u)-\Phi(v)}{d(\hat u,v)}.
\]
Therefore there exists a neighborhood $\mathcal V$ of $\hat u$ such that
\[
    \sup_{v\in\mathcal V\setminus\{\hat u\}} \frac{\Phi(\hat u)-\Phi(v)}{d(\hat u,v)} < \frac{\varepsilon}{\lambda}.
\]
From this, for all $v\in\mathcal V\setminus\{\hat u\}$, 
\begin{align*}
    \Phi(\hat u) = \Phi(v) + \frac{\Phi(\hat u)-\Phi(v)}{d(\hat u,v)} d(\hat u,v) < \Phi(v) + \frac{\varepsilon}{\lambda} d(\hat u,v). 
\end{align*}
This completes the proof. \hfill\(\square\) 



\subsubsection{Proof of Corollary \ref{cor_gms}}
We first note that the lower semicontinuity of $\Phi$ guarantees that the global metric slope 
$|\partial\Phi|$ is a strong upper gradient of $\Phi$  \cite[Theorem~1.2.5]{AGS_2008} (the result does not require completeness of $X$). Moreover, the argument in the proof of \cite[Proposition~2.2]{LT_2024} shows that 
$|\partial\Phi|:X\to\mathbb R\cup\{+\infty\}$ is lower semicontinuous (although the result is stated there for real-valued functions, the same argument 
applies to extended real-valued functions).  Now, by Theorem \ref{thm_gvp_alt}, one can find $\hat u\in\mathcal M(|\partial\Phi|)$ such that $d(u,\hat u)\le\lambda$, $\Phi(\hat u)+\frac\varepsilon\lambda d(u,\hat u)\le\Phi(u)$ and $|\partial\Phi|(\hat u)\le\varepsilon/\lambda$. 
By definition of metric slope, 
\[
    \frac\varepsilon\lambda\ge|\partial\Phi|(\hat u) \ge  \sup_{v\in X\setminus\{\hat u\}}\frac{\Phi(\hat u)-\Phi(v)}{d(\hat u,v)}.
\]
From this, for all $v\in X\setminus\{\hat u\}$, 
\begin{align*}
    \Phi(\hat u) = \Phi(v) + \frac{\Phi(\hat u)-\Phi(v)}{d(\hat u,v)} d(\hat u,v) \le \Phi(v) + \frac{\varepsilon}{\lambda} d(\hat u,v). 
\end{align*}
The result follows. \hfill\(\square\)



\section{Some standard consequences of the variational principle}

\subsection{A weak Palais--Smale condition}
A standard use of variational principles is to derive the existence of minimizers under Palais--Smale type conditions \cite{KP_2000, MW_2010, PS_1964, S_2008}. The mechanism is simple; one produces a sequence with nearly minimal energy and vanishing gradient, and then extracts a critical point. 
\smallbreak\noindent
Let $(X, d)$ be a metric space and   $\Phi:X\to\mathbb R\cup\{+\infty\}$ be a proper functional bounded from below. 
\smallbreak \noindent
We say that an upper gradient  $g:X\to\mathbb R\cup\{+\infty\}$ of $\Phi$ satisfies the \textit{weak Palais--Smale condition} if for every sequence $(u_n)_{n\in\mathbb N}\subseteq X$ satisfying
    \begin{align*}
        \sup_{n\in\mathbb N}\Phi(u_n)<+\infty\qquad\text{and}\qquad \lim_{n\longrightarrow+\infty}g(u_n) = 0
    \end{align*}
    there exists $u\in X$ with  $g(u) = 0$ such that 
    \[
    \displaystyle\liminf_{n\longrightarrow+\infty}\Phi(u_n)\le\Phi(u)\le\displaystyle\limsup_{n\longrightarrow+\infty}\Phi(u_n).
    \]
This version of the Palais--Smale condition is enough to conclude the existence of minimizers when the functional admits enough gradient flows. 
\begin{proposition}\label{P1}
    Let  $g:X\to\mathbb R\cup\{+\infty\}$ be an upper gradient  of $\Phi$  such that  $[\Phi\le\alpha]\subseteq\mathcal M(g)$ for some $\alpha>\inf\Phi$.
  If $g$ satisfies the weak Palais--Smale condition, then $\argmin_{v\in X}\Phi(v)\neq\emptyset$. 
\end{proposition}
\noindent
We now make some remarks on this proposition.
\begin{remark}
No completeness of the underlying space is assumed; the compactness requirement is encoded in the weak Palais--Smale condition for the upper gradient.
\end{remark}

\begin{remark}
We mention that in \cite[Theorem 3.4]{CDW_2022} a gradient flow--type argument is given to prove the existence of minimizers under a generalized Palais--Smale condition. Such a result can be derived directly from the variational principles developed in the previous section.
\end{remark}



\subsection{A metric deformation lemma}
A standard geometric use of variational principles is to build deformations of sublevel sets when no critical points are present \cite{C_2013,CDM_1993,P_1986}. In the smooth setting, if the size of the gradient is bounded from below by a positive constant on an intermediate range of energy values, then the negative gradient flow decreases the energy at a uniform rate. 
\smallbreak\noindent
During this subsection, $(X,d)$ denotes a metric space and $\Phi:X\to\mathbb R\cup\{+\infty\}$ a proper functional bounded from below.  
Fix $\alpha,\beta\in(\inf\Phi,+\infty)$ such that $\alpha<\beta$. 
\smallbreak\noindent
We have the following metric
deformation lemma.
\begin{proposition}\label{P2}
Let $g:X\to\mathbb R\cup\{+\infty\}$ be an upper gradient of $\Phi$ such that $[\alpha<\Phi\le\beta]\subseteq\mathcal M(g)$.  
  Assume that there exists $\sigma>0$ such that 
\begin{align}\label{hdl}
    g(v)\ge\sigma \quad\forall v\in[\alpha\le\Phi\le\beta].
\end{align}
Then, for every $u\in [\Phi\le\beta]$   there exists an absolutely continuous curve $\eta:[0,1]\to X$ such that 
	\begin{align*}
		(a)\,\,\, \eta(0)=u\qquad\,\,\,\, (b)\,\,\,\eta(1)  \in [\Phi\le\alpha]\qquad \,\,\,\,(c)\,\,\,   d(u,\eta(t)) \le \frac{t}{\sigma}\max\{0,\Phi(u)-\alpha\}\quad \forall t\in[0,1]. 
	\end{align*}
\end{proposition}
\noindent
We make the following remark. 
\begin{remark}
It follows from the proof of Proposition \ref{P2} that every $u\in[\alpha<\Phi\le\beta]$ can be joined to the sublevel set $[\Phi\le\alpha]$ by a gradient-flow path. Under additional hypotheses ensuring continuous dependence of gradient flows on the initial datum, and extending this construction by the identity on $[\Phi\le\alpha]$, one could obtain a deformation retract of $[\Phi\le\beta]$ onto $[\Phi\le\alpha]$.
\end{remark}

\subsection{A local error bound for critical points}
In this subsection we show how quantitative control of the set of critical points in terms of the metric slope can yield a local growth estimate.
\smallbreak\noindent
Throughout the subsection, we consider  a metric space  $(X, d)$ and  a proper functional  $\Phi:X\to\mathbb R\cup\{+\infty\}$ bounded from below. We assume that the local metric slope $|\nabla\Phi|:X\to\mathbb{R}\cup\{+\infty\}$ is an upper gradient of $\Phi$. The set of critical points of $\Phi$ is defined by
\[
    \operatorname{crit}\Phi:=\{u\in X:\ |\nabla\Phi|(u)=0\}.
\]
We fix $\theta>0$ and a reference minimizer $u^*\in X$ of $\Phi$; observe that $u^*\in\operatorname{crit}\Phi$.

\begin{proposition}\label{Pcrit}
 Suppose that there exist a neighborhood $\mathcal V$ of $u^*$ and $\kappa>0$ such that
\begin{align}\label{leb}
    \operatorname{dist}\bigl(u,\operatorname{crit}\Phi\bigr)\le \kappa\,|\nabla\Phi|(u)^{\theta}
    \quad\forall u\in \mathcal V\cap\mathcal M(|\nabla\Phi|).
\end{align}
Then there exists a neighborhood $\mathcal W$ of $u^*$ and $c>0$ such that
\begin{align}\label{et}
    \Phi(u)-\Phi(u^*) \;\ge\; c\,\operatorname{dist}\bigl(u,\operatorname{crit}\Phi\bigr)^{1+\frac{1}{\theta}}
    \quad\forall u\in \mathcal W\cap\mathcal M(|\nabla\Phi|).
\end{align}
\end{proposition}

\noindent
We now make some remarks. 
\begin{remark}
The estimate (\ref{leb}) 
can be regarded as a local \emph{error bound}; the slope acts as a residual, and small residual forces
proximity to the critical set. In a normed space, when the functional \(\Phi\) is continuously differentiable, this reads
\(\operatorname{dist}(u,\operatorname{crit}\Phi)\lesssim \|D\Phi(u)\|^{\theta}\), which is a form of
Hölder metric subregularity for the gradient mapping; see, e.g., \cite{A_2008,AC_2014,DR_2014,JOV_2025,MN_2013}.
\end{remark}

\begin{remark}
    The growth condition (\ref{et}) is a well-known error bound; see \cite{AC_2014,AC_2017,K_2015}. In  Wasserstein spaces, it is common to refer to conditions such as (\ref{et}) as entropy–transport inequalities or Talagrand inequalities; see \cite[Section~3.8]{HM_2019} and \cite[Section~3.1]{BB_2018}.
\end{remark}

\begin{remark}
The proof yields an explicit constant, in the sense that for any $c\in(0,c_0)$ there exists a neighborhood $\mathcal W$ so that (\ref{et}) holds, where
\[
c_0=\frac{\theta}{(1+\theta)^{1+1/\theta}}\,\kappa^{-1/\theta}.
\]
In particular, the dependence on \(\kappa\) and \(\theta\) is completely explicit.
\end{remark}


\subsection{A Caristi fixed point theorem}
The fixed point theorem of Caristi \cite{C_1976} is a classical consequence of the variational principle, see, e.g., the monograph \cite{KS_2014}. It asserts that if a mapping always decreases the energy by at least the distance it moves the point, then it has a fixed point. In this subsection we obtain an analogous fixed point result based on Theorem~\ref{thm_gvp}.
\smallbreak\noindent
Let $(X, d)$ be a metric space and   $\Phi:X\to\mathbb R\cup\{+\infty\}$ be a proper functional bounded from below. Let us fix $\gamma\ge0$.  Consider $|\partial\Phi|_\gamma:X\to\mathbb R\cup\{+\infty\}$ given by 
\begin{align}\label{alphsup}
    |\partial\Phi|_\gamma(u) :=\left\{
    \begin{array}{lcc}
    \max\Biggl\{
\displaystyle\sup_{v\in X\setminus\{u\}}
\biggl(
\frac{\Phi(u)-\Phi(v)}{d(u,v)} - \frac{\gamma}{2}\,d(u,v)
\biggr),
\,0
\Biggr\} & \text{if} & u\in\dom\Phi \\ \\
    +\infty & \text{if} & u\notin\dom\Phi.
    \end{array}
    \right.
\end{align}
Notice that $|\partial\Phi|_\gamma$ is nonnegative and reduces to the  global slope when $\gamma=0$. 

\begin{proposition}\label{P3}
   Suppose that   $|\partial\Phi|_\gamma:X\to\mathbb R\cup\{+\infty\}$ is a strong upper gradient  of $\Phi$  and that $[\Phi<\inf \Phi+2/\gamma]\cap\mathcal M(|\partial\Phi|_\gamma)\neq\emptyset$.
    If $T: X \to X$ is a mapping satisfying
    \begin{align}\label{cariscon}
            \Phi(v) \ge \Phi\big(T(v)\big) + d\big(v, T(v)\big) \quad \forall v\in X,
    \end{align}
then there exists $\hat u\in\mathcal M(|\partial\Phi|_\gamma)$ such that $T(\hat u) = \hat u$.
\end{proposition}
\noindent
We now make some remarks. 

\begin{remark}
    It is expected that $|\partial\Phi|_\gamma$ be a strong upper gradient under geodesic semi-convexity assumptions \cite[Proposition 2.22]{HM_2019}. In \cite[Definition 2.7]{MS_2020}, (\ref{alphsup}) is referred to as the \textit{global $\gamma$-slope}. 
\end{remark}

\begin{remark}
In the case \(\gamma=0\) we interpret the threshold \(2/\gamma\) as \(+\infty\); hence
the assumption \( [\Phi<\inf\Phi+2/\gamma]\cap\mathcal M(|\partial\Phi|_\gamma)\neq\emptyset \)
simply reduces to \(\mathcal M(|\partial\Phi|)\neq\emptyset\).
\end{remark}
\begin{remark}
With verbatim adaptations, more general strong upper gradients can be considered in  Proposition \ref{P3}, replacing the discounted linear term by a modulus of continuity vanishing at zero. These upper gradients arise in the study of nonconvex Wasserstein gradient flows; see \cite[Propositions 2.6 and 2.7]{C_2017}.
\end{remark}

\subsection{Proofs}

\subsubsection{Proof of Proposition \ref{P1}}
    Let $(u_n)_{n\in\mathbb N}$ be a sequence so that $\Phi(u_n)\le \inf\Phi+1/n$ for all $n\in\mathbb N$. By Theorem \ref{thm_gvp}, there exists a sequence $(\hat u_n)_{n\in\mathbb N}$ such that $\Phi(\hat u_n)\le \Phi(u_n)\le \alpha$ and $g(\hat u_n)\le 1/n$  for all $n\in\mathbb N$ large enough. From this and the assumption that $g$ satisfies the weak Palais--Smale condition, we conclude the existence of $\hat u\in X$ satisfying 
    \[
        \Phi(\hat u)\le \limsup_{ n\longrightarrow+\infty}\Phi(\hat u_n)\le \limsup_{ n\longrightarrow+\infty}\Phi(u_n) =\inf\Phi.
    \]
    Whence the result follows. \hfill\(\square\) 
    
\subsubsection{Proof of Proposition \ref{P2}}

Let $u\in [\Phi\le\beta]$ be given. If $u\in[\Phi\le\alpha]$, then $\eta:[0,1]\to X$ given by $\eta(t)=u$ satisfies all the claimed properties. Suppose that $u\in[\alpha<\Phi\le\beta]$. Since $u\in[\alpha<\Phi\le\beta]\subseteq\mathcal M(g)$, there exists  a gradient flow
$v:(0,+\infty)\to X$ of $\Phi$ with $v(t)\longrightarrow u$ and $\Phi(v(t))\longrightarrow\Phi(u)$ as $t\longrightarrow 0^+$.
Set $\tau:=\inf\{t>0:\,\Phi(v(t))\le\alpha\}$. Observe that $\Phi(v(\tau))=\alpha$. 
By Lemma~\ref{AB},
\[
\Phi(u)-\Phi(v(\tau))=\int_0^{\tau} g(v(s))^2\,ds\ge \tau \sigma^2.
\]
From this, 
\(\tau\le(\Phi(u)-\alpha)/\sigma^2\).
Define $\eta:[0,1]\to X$ by $\eta(t):=v(\tau t^2)$. Then $\eta$ is an absolutely continuous curve satisfying $\eta(0)=u$  and $\eta(1)\in[\Phi\le\alpha]$.
Since $\Phi(\eta(t))\ge\alpha$ for all $t\in[0,1]$, we get 
\begin{align*}
    d(u,\eta(t)) \;&\le\; \int_0^{\tau t^2} |\dot v|(r)\,dr
\;\le\; \sqrt{\tau t^2}\left(\int_0^{\tau t^2} |\dot v|^2(r)\,dr\right)^{1/2}
\;=\; t\sqrt{\tau}\,\big(\Phi(u)-\Phi(\eta(t))\big)^{1/2}\\ 
&\;\le t\sqrt{\tau}\big(\Phi(u)-\alpha\big)^{1/2}\le  t \frac{\Phi(u)-\alpha}{\sigma}
\end{align*}
for all $t\in[0,1]$;  whence the result follows. \hfill\(\square\) 


\subsubsection{Proof of Proposition \ref{Pcrit}}
Set $\beta:=\theta(1+\theta)^{-1}$ and 
\[
c_0:=\left(\frac{\beta^\theta(1-\beta)}{\kappa}\right)^{1/\theta}
= \frac{\theta}{(1+\theta)^{1+1/\theta}}\,\kappa^{-1/\theta}.
\]
Fix $\delta\in(0,1)$. We prove that there exists a neighborhood $\mathcal W$ of $u^*$ such that
\begin{equation}\label{eq:goal}
    \Phi(u)-\Phi(u^*)
    \;\ge\; (1-\delta)c_0\operatorname{dist}\bigl(u,\operatorname{crit}\Phi\bigr)^{1+1/\theta}
    \quad\forall u\in \mathcal W\cap \mathcal M(|\nabla\Phi|).
\end{equation}
The desired conclusion then will follow with $c=(1-\delta)c_0$. 
Suppose  that \eqref{eq:goal} fails. Then we can find a sequence $(u_n)_{n\in\mathbb N}\subseteq\mathcal M(|\nabla\Phi|)$ such that
\[
u_n\longrightarrow u^*
\quad\text{and}\quad
\Phi(u_n)-\Phi(u^*)<(1-\delta)c_0\operatorname{dist}\bigl(u_n,\operatorname{crit}\Phi\bigr)^{1+1/\theta}
\quad\forall n\in\mathbb N.
\]
Observe that, for  any $n\in\mathbb N$, $u_n\notin\overline{\operatorname{crit}\Phi}$; otherwise $\Phi(u_n)<\Phi(u^*)=\inf\Phi$. 
 For each $n\in\mathbb N$, set 
\[
\varepsilon_n:=c_0(1-\delta)\operatorname{dist}(u_n,\operatorname{crit}\Phi)^{1+1/\theta}
\quad\text{and}\quad 
\lambda_n:=\beta \operatorname{dist}(u_n,\operatorname{crit}\Phi).
\]
We can apply Theorem~\ref{thm_gvp} to obtain a sequence $(\hat u_n)_{n\in\mathbb N}$ such that 
\begin{equation}\label{eq:vp-out}
    d(u_n,\hat u_n)\le\lambda_n\quad\text{and}\quad
    |\nabla\Phi|(\hat u_n)\le\frac{\varepsilon_n}{\lambda_n}\quad \forall n\in\mathbb N. 
\end{equation}
This implies that 
\begin{equation*}
    \hat u_n\longrightarrow u^*\quad\text{and}\quad|\nabla\Phi|(\hat u_n)\le \frac{(1-\delta)c_0}{\beta}\,\operatorname{dist}(u_n,\operatorname{crit}\Phi)^{1/\theta}\quad\forall n\in\mathbb N
\end{equation*}
For $n\in\mathbb N$ large enough,  $\hat u_n\in\mathcal V$ and hence
\begin{equation*}
   \operatorname{dist}(\hat u_n,\operatorname{crit}\Phi)\le\kappa\Big(\frac{(1-\delta)c_0}{\beta}\Big)^\theta \operatorname{dist}(u_n,\operatorname{crit}\Phi).
\end{equation*}
On the other hand, by the triangle inequality and \eqref{eq:vp-out},
\begin{align*}
\operatorname{dist}\bigl(u_n,\operatorname{crit}\Phi\bigr)
&\le d(u_n,\hat u_n)+\operatorname{dist}\bigl(\hat u_n,\operatorname{crit}\Phi\bigr)
\le \lambda_n + \operatorname{dist}\bigl(\hat u_n,\operatorname{crit}\Phi\bigr)\\
&\le \beta\,\operatorname{dist}\bigl(u_n,\operatorname{crit}\Phi\bigr)
     + \kappa\Big(\frac{(1-\delta)c_0}{\beta}\Big)^\theta
       \operatorname{dist}\bigl(u_n,\operatorname{crit}\Phi\bigr)
\end{align*}
for all $n\in\mathbb N$ large enough. This implies
\[
1 \le \beta + \kappa\Big(\frac{(1-\delta)c_0}{\beta}\Big)^\theta = \beta +  (1-\delta)^\theta(1-\beta)<\beta + (1-\beta)=1.
\]
This is the sought contradiction. \hfill\(\square\)


\subsubsection{Proof of Proposition \ref{P3}}
By assumption, there exists $u\in \mathcal M(|\partial\Phi|_\gamma)$ satisfying $\Phi(u)<\inf\Phi+2/\gamma$. If $\Phi(u)=\inf\Phi$, the result follows trivially as
\[
    d(u,Tu)\le \Phi(u)-\Phi(Tu)=\inf\Phi-\Phi(Tu)\le 0. 
\]
Assume then that $\Phi(u)\neq\inf\Phi$. Set $\varepsilon:=\Phi(u)-\inf\Phi$ and note that $\varepsilon\in(0,2/\gamma)$. Choose $\lambda>0$ so that
\begin{align}\label{lamch}
    \frac{\varepsilon}{\lambda}<1-\frac{\gamma}{2}\,\varepsilon.
\end{align}
By Theorem~\ref{thm_gvp}, there exists $\hat u\in\mathcal M(|\partial\Phi|_\gamma)$ such that $\Phi(\hat u)\le\Phi(u)$ and 
\[
    \Phi(\hat u)\le \Phi(v)+\frac{\varepsilon}{\lambda}\,d(v,\hat u)+\frac{\gamma}{2}\,d(v,\hat u)^2 \quad\forall v\in X\setminus\{\hat u\}.
\]
Combining this with \eqref{cariscon}, we obtain
\[
d(\hat u,T\hat u)\le \Phi(\hat u)-\Phi(T\hat u)
\le \frac{\varepsilon}{\lambda}\,d(\hat u,T\hat u)+\frac{\gamma}{2}\,d(\hat u,T\hat u)^2.
\]
Hence $\big(1-\varepsilon\lambda^{-1}\big)d(\hat u,T\hat u)\le \gamma/2\, d(\hat u, T\hat u)^2$. If $d(\hat u, T\hat u)>0$, then
\[
    1-\frac{\varepsilon}{\lambda} \le \frac{\gamma}{2} d(\hat u,T\hat u)\le\frac{\gamma}{2}\big( \Phi(\hat u)-\Phi(T\hat u)\big)\le \frac{\gamma}{2}\big(\Phi(\hat u)-\inf\Phi\big)\le\frac{\gamma}{2} \varepsilon.
\]
This is a contradiction to (\ref{lamch}); whence it follows that $d(\hat u,T\hat u)=0$. \hfill\(\square\)






\section{Mabuchi's energy and the Calabi flow}
\noindent
In this section we illustrate one of the variational principles (Corollary~\ref{cor_gms}) established in
Section~\ref{sec:mgf} for the Mabuchi $K$-energy on the space of K\"ahler potentials.
\smallbreak\noindent
Throughout this section, we fix a compact K\"ahler manifold $(M^{2n},\omega,J)$.

\subsection{The Calabi flow}
The space of K\"ahler potentials in the cohomology class $[\omega]$ is given by
\[
\mathcal H:=\Bigl\{\phi\in C^\infty(M)\;\Big|\;\omega+\sqrt{-1}\,\partial\bar\partial\phi>0\Bigr\}.
\]
 For \(\phi\in\mathcal H\),  we denote by
\(\omega_\phi:=\omega+\sqrt{-1}\,\partial\bar\partial\phi\) the induced K\"ahler form, and by $s_\phi$ the scalar curvature of the K\"ahler metric associated to $\omega_\phi$.  We set 
\[
\bar s:=\frac{1}{\operatorname{vol}(\omega_\phi)}\int_M s_\phi\,\omega_\phi^n,
\]
and note that \(\bar s\) is independent of the K\"ahler potential \(\phi\in\mathcal H\).
\smallbreak\noindent
We say that a one-parameter family $(\phi_t)_{t\in[0,+\infty)}\subset\mathcal H$ is a \textit{smooth solution of the Calabi flow} if the mapping  $[0,+\infty)\times M\ni(t,z)\longmapsto \phi_t(z)\in\mathbb R$ is smooth and 
\begin{align}\label{Calabi}
     \frac{\partial}{\partial t}\phi_t = s_{\phi_t} - \bar s\quad\,\,\text{in}\,\, (0,+\infty)\times M.
\end{align}
For such a solution, we call $\phi_0\in\mathcal H$ the \textit{initial datum}. It is still an open conjecture whether solutions exist in general. However, their existence has been established in a number of settings, e.g., 
\begin{itemize}[label=\raisebox{0.2ex}{\scriptsize\(\bullet\)}]
  \item (\textit{Complex dimension one}). On a compact Riemann surface, the Calabi flow exists for all time and converges to a constant scalar curvature metric; see Chru\'sciel \cite{C_1991B} and Chen \cite{C_2001}. 

  \item (\textit{Stability near cscK}). If $[\omega]$ contains a cscK metric, then for initial potentials sufficiently small in $C^{3,\alpha}$ (with respect to the cscK metric), the Calabi flow exists for all time and converges in $C^\infty$ exponentially fast to a cscK metric in $[\omega]$ (unique up to automorphisms) \cite[Theorem 1.2]{CH_2008}. 

  \item (\textit{Toric Fano surfaces}). On a toric Fano surface with positive extremal Hamiltonian potential, starting from a toric invariant initial metric whose Calabi energy satisfies an explicit smallness condition, it was proved that the Calabi flow exists for all time and converges (along a subsequence) to an extremal metric in the Cheeger--Gromov sense \cite[Theorem 1.1]{CH_2010}. 
  
  \item (\textit{Flat complex tori}). Consider the space  $\mathbb C^n/(\mathbb Z^n+i\mathbb Z^n)$ and restrict to $T^n$-invariant K\"ahler metrics in a fixed class. It was proved that in complex dimension $n=2$, for any initial datum in their invariant class, the Calabi flow exists for all time \cite[Theorem 1.2]{FH_2012}.
\end{itemize}

\subsection{Mabuchi's energy}
We begin by recalling the weak Riemannian structure on the space of K\"ahler potentials  \cite{D_1999,M_1987,S_1992}. 
Given $\phi\in\mathcal H$, the tangent space $T_\phi\mathcal H$ can be identified with $C^\infty(M)$. 
\smallbreak\noindent
For $\phi\in\mathcal H$ and $\psi_1,\psi_2\in T_\phi\mathcal H$,
consider the inner product on $T_\phi\mathcal H$ given by
\begin{align}\label{eq:mabuchi-metric}
\langle \psi_1,\psi_2\rangle_\phi:=\int_M \psi_1\psi_2\,\omega_\phi^n.
\end{align}
This induces a weak Riemannian metric, and hence a (path-length) distance $d$ on $\mathcal H$. It is known that $d$ is nondegenerate, so $(\mathcal H,d)$ is a metric space \cite[Theorem~6]{C_2000}; see also \cite[Theorem~1.2]{B_2012}.
\smallbreak\noindent
The \textit{K-energy} is the functional  $\nu:\mathcal H\to\mathbb R$ defined by
\begin{align}\label{eq:k-energy-def}
\nu(\phi):=-\int_0^1\int_M \dot\alpha_t\,\bigl(s_{\alpha_t}-\bar s\bigr)\,\omega_{\alpha_t}^n\,dt,
\end{align}
where $(\alpha_t)_{t\in[0,1]}\subset\mathcal H$ is any smooth one-parameter family with $\alpha_0=0$ and
$\alpha_1=\phi$. It is a standard fact that the right-hand side of \eqref{eq:k-energy-def} is independent of the
choice of $(\alpha_t)_{t\in[0,1]}$. 
\smallbreak\noindent
The Calabi flow is expected to be the gradient flow of the K-energy; in the following lemmas, we make this precise in the metric gradient flow framework of Subsection \ref{mgfs}.
\begin{lemma}\label{lem:slope-Kenergy}
Let $g:\mathcal H\to\mathbb R$ be given by 
\begin{align}\label{sugc}
    g(\phi):=
\left(\int_M \bigl(s_\phi-\bar s\bigr)^2\,\omega_\phi^n\right)^{1/2}.
\end{align}
Then, $g$ is a strong upper gradient of $\nu$. Moreover, $|\nabla\nu|(\phi)=|\partial\nu|(\phi)=g(\phi)$ for all $\phi\in\mathcal H$.
\end{lemma}

\noindent
Before passing to the next lemma, we make some remarks.  
\begin{remark}\label{rcelsc}
    It is known that the $K$-energy is lower semicontinuous \cite[Lemma 5.11]{S_2014}, and therefore by \cite[Proposition 2.2]{LT_2024} its global metric slope is also lower semicontinuous. In particular, by Lemma \ref{lem:slope-Kenergy}, the upper gradient $g:\mathcal H\to\mathbb R$ in (\ref{sugc}) is lower semicontinuous. 
\end{remark}

\begin{remark}
The property \( |\nabla\nu|=|\partial\nu| \) has been previously referred to as the \emph{slope cone property} in the literature
\cite[Definition~3.1]{CG_2012}. We point out that the \(K\)-energy is not geodesically convex in the usual sense
\cite[Definition~2.4.3]{AGS_2008}, though it can be regarded as convex in weaker notions; see \cite{BB_2017,BD_2017}.
\end{remark}

\begin{remark}
Let $\phi\in\mathcal H$. The quantity
\begin{align}\label{calener}
    \mathcal C(\phi):=\int_M (s_\phi-\bar s)^2\,\omega_\phi^n
\end{align}
is the so-called Calabi energy of $\omega_\phi$. It follows from (\ref{sugc}) and (\ref{calener}) that  $g(\phi)=\mathcal C(\phi)^{1/2}$; in particular, $g(\phi)=0$ if and only if $\omega_\phi$ has constant scalar curvature.
\end{remark}

\begin{lemma}\label{lem:gfc}
    Let $(\phi_t)_{t\in[0,+\infty)}$ be a smooth solution of the Calabi flow (\ref{Calabi}). Then, the curve $(0,+\infty)\ni t\longmapsto\phi_t\in\mathcal H$ is a gradient flow of the $K$-energy with respect to the upper gradient $g:\mathcal H\to\mathbb R$ in (\ref{sugc}) in the sense of Definition \ref{dgf}.
\end{lemma}

\subsection{The variational principle}

We are now ready to state the main result of this section, which provides a substitute for the classical Ekeland variational principle obtained from the Calabi flow.

\begin{theorem}\label{thmcal}
Assume that the $K$-energy (\ref{eq:k-energy-def}) is bounded from below, and that there exists a smooth solution
$(\phi_t)_{t\in[0,+\infty)}$ of the Calabi flow. 
Let $\varepsilon>0$ be such that
\[
\nu(\phi_0)\le \inf_{\psi\in\mathcal H}\nu(\psi)+\varepsilon.
\]
Then for any $\lambda>0$ there exists $\tau\in[0,+\infty)$ such that
\begin{align*}
(a)\;\; d(\phi_0,\phi_\tau)\le\lambda\qquad
(b)\;\; \nu(\phi_\tau)+\frac{\varepsilon}{\lambda}\,d(\phi_0,\phi_\tau)\le \nu(\phi_0)\qquad
(c)\;\; \nu(\phi_\tau)\le \nu(\psi)+\frac{\varepsilon}{\lambda}\,d(\psi,\phi_\tau)\quad\forall \psi\in\mathcal H.
\end{align*}
\end{theorem}

\noindent
We close the section with a few remarks on the theorem. 

\begin{remark}
If the class \([\omega]\) contains a cscK metric, then the \(K\)-energy is automatically bounded from below.
Intuitively, boundedness from below is expected to reflect a stability property; it was shown in
\cite[Theorem~4.1]{DR_2017} that it implies \(K\)-semistability \cite[Definition~2.22]{DR_2017} of the K\"ahler manifold.
\end{remark}

\begin{remark}
Note that item $(c)$ in Theorem~\ref{thmcal} is equivalent to the estimate 
\(\mathcal C(\phi_\tau)\le \varepsilon^2/\lambda^2\), where \(\mathcal C(\phi_\tau)\) denotes the Calabi energy
\eqref{calener} of the K\"ahler form \(\omega_{\phi_\tau}\).
\end{remark}

\begin{remark}
The classical variational principle can be applied to the metric completion of $(\mathcal H,d)$, the so-called Darvas space \cite[Theorem~1]{D_2017}. This is because the $K$-energy admits a lower semicontinuous extension to this completion \cite[Definition~4.2 and Lemmas~5.10--5.11]{S_2014}. Therefore, the classical Ekeland principle can give a point in the completion satisfying items $(a)$--$(c)$ of Theorem~\ref{thmcal}. The grace of Theorem~\ref{thmcal} is that it additionally provides a point in the smooth space $\mathcal H$, selected along the Calabi-flow trajectory.
\end{remark}

\subsection{Proofs}
\subsubsection{Proof of Lemma \ref{lem:slope-Kenergy}}
By definition of the local and global metric slopes,
\begin{align}\label{natin}
    |\nabla\nu|(\phi)\le|\partial\nu|(\phi)\qquad\forall\,\phi\in\mathcal H.
\end{align}
We prove below that \( |\partial\nu|(\phi)\le g(\phi)\) and \(|\nabla\nu|(\phi)\ge g(\phi)\) for all \(\phi\in\mathcal H\); combining these inequalities with \eqref{natin}, we can then conclude the desired identity. The fact that \(g:\mathcal H\to\mathbb R\) is a strong upper gradient is then a consequence of \cite[Theorem 1.2.5]{AGS_2008} (the argument in the proof of \cite[Theorem 1.2.5]{AGS_2008} does not require completeness of the underlying space) and that the $K$-energy is lower semicontinuous \cite[Lemma 5.11]{S_2014}.
\smallbreak\noindent
\textbf{(Inequality $|\partial\nu|\le g$}). For any $\phi\in\mathcal H$, 
 by \cite[Theorem~1.2]{C_2009}, 
\begin{equation*}
\nu(\psi)\ge\nu(\phi)-d(\phi,\psi)\,g(\phi)\quad\forall\psi\in\mathcal H.
\end{equation*}
Rearranging and taking supremum, we obtain that 
\[
   \sup_{\psi\in \mathcal H\setminus\{\phi\}} \frac{\nu(\phi)-\nu(\psi)}{d(\phi,\psi)}\le g(\phi)\quad\forall\phi\in\mathcal H.
\]
Taking into account that $g:\mathcal H\to\mathbb R$ is nonnegative, we get $|\partial\nu|(\phi)\le g(\phi)$ for all $\phi\in\mathcal H$.
\smallbreak\noindent
\textbf{(Inequality $|\nabla\nu|\ge g$}).
We argue along the lines of the argument in \cite[Lemma~3.2]{J_2016}.  Fix $\phi\in\mathcal H$ arbitrary. 
If $g(\phi)=0$, the inequality holds automatically; we may therefore assume
$g(\phi)>0$.
\smallbreak\noindent
By standard results on short-time existence for the Calabi flow (see, e.g., \cite[Section~3]{CH_2008}), there exists
$T>0$ and a smooth one-parameter family $(\phi_t)_{t\in[0,T)}\subset\mathcal H$ such that $\phi_0=\phi$ and 
\[
 \frac{\partial}{\partial t}\phi_t = s_{\phi_t} - \bar s\quad\,\,\text{in}\,\, (0,T)\times M.
\]
From the first variation formula for the $K$-energy, 
\[
\frac{d}{dt}\nu(\phi_t)
=-\int_M \dot\phi_t\,(s_{\phi_t}-\bar s)\,\omega_{\phi_t}^n=-\int_M (s_{\phi_t}-\bar s)^2\,\omega_{\phi_t}^n
=-g(\phi_t)^2
\]
for all $t\in[0,T)$. 
Integrating then yields
\begin{equation}\label{eq:nu-drop}
\nu(\phi_0)-\nu(\phi_t)=\int_0^t g(\phi_s)^2\,ds
\qquad\forall t\in[0,T).
\end{equation}
On the other hand, by definition of the distance $d$ on $\mathcal H$,
\begin{equation}\label{eq:d-bound}
d(\phi_0,\phi_t)\le \int_0^t \biggl(\int_M \dot\phi_s^{\,2}\,\omega_{\phi_s}^n\biggr)^{1/2}\,ds
=\int_0^t g(\phi_s)\,ds
\qquad\forall t\in[0,T).
\end{equation}
Combining \eqref{eq:nu-drop} and \eqref{eq:d-bound} we obtain
\[
\frac{\nu(\phi_0)-\nu(\phi_t)}{d(\phi_0,\phi_t)}
\;\ge\;
\frac{\displaystyle\int_0^t g(\phi_s)^2\,ds}
     {\displaystyle\int_0^t g(\phi_s)\,ds}\quad \forall t\in(0,T).
\]
Since the scalar curvature depends smoothly on the metric, 
\[
\lim_{t\longrightarrow0^+}\frac{1}{t}\int_0^t g(\phi_s)^2\,ds=g(\phi_0)^2\quad\text{and}\quad
\lim_{t\longrightarrow0^+}\frac{1}{t}\int_0^t g(\phi_s)\,ds=g(\phi_0).
\]
We then conclude that 
\[
\liminf_{t\longrightarrow0^+}\frac{\nu(\phi_0)-\nu(\phi_t)}{d(\phi_0,\phi_t)}
\;\ge\;
\frac{g(\phi_0)^2}{g(\phi_0)}
=g(\phi_0).
\]
It follows that
\[  |\nabla\nu|(\phi_0)\ge\limsup_{\psi\longrightarrow\phi_0}\frac{\nu(\phi_0)-\nu(\psi)}{d(\phi_0,\psi)}\ge \liminf_{t\longrightarrow0^+}\frac{\nu(\phi_0)-\nu(\phi_t)}{d(\phi_0,\phi_t)}
\;\ge\;
\frac{g(\phi_0)^2}{g(\phi_0)}
=g(\phi_0).
\]
Since $\phi_0=\phi$ and $\phi\in\mathcal H$ was arbitrary, this concludes the proof. \hfill\(\square\)

\subsubsection{Proof of Lemma \ref{lem:gfc}}
    We begin by observing that, by Lemma \ref{lem:slope-Kenergy}, the function $g:\mathcal H\to\mathbb R$ defined in \eqref{sugc} is indeed a strong upper
gradient of $\nu
:\mathcal H\to\mathbb R$. Let $(\phi_t)_{t\in[0,+\infty)}$ be a smooth solution of the Calabi flow  \eqref{Calabi}. Consider the function $\alpha:(0,+\infty)\to\mathcal H$ given by $\alpha(t):=\phi_t$.  Since $\dom\nu=\mathcal H$, the function $\alpha$ satisfies item $(1)$ of Definition \ref{dgf}; below we prove that item $(2)$ also holds. 
\smallbreak\noindent
Since $d$ is the path-length distance induced by the metric 
\eqref{eq:mabuchi-metric}, 
\begin{align}\label{absconal}
    d(\phi_s,\phi_t)
\le \int_s^t \Bigl(\int_M \dot\phi_r^{\,2}\,\omega_{\phi_r}^n\Bigr)^{1/2}\,dr=\int_s^t\Bigl(\int_M (s_{\phi_r}-\bar s)^2\,\omega_{\phi_r}^n\Bigr)^{1/2}\, dr
=\int_s^tg(\phi_r)\,dr
\end{align}
for all $s,t\in[0,+\infty)$ satisfying $s\le t$.  Therefore the curve 
 $\alpha:(0,+\infty)\to\mathcal H$ is locally absolutely continuous in $(\mathcal H,d)$. It also follows from (\ref{absconal}) that 
its metric  derivative exists almost everywhere, belongs to $L^1_{\operatorname{loc}}(0,+\infty)$, and satisfies
\[
    |\dot\alpha|(t)\le g(\phi_t)\quad\text{for almost every $t\in(0,+\infty)$.}
\]
The first variation formula for the $K$-energy yields
\begin{align}\label{eq:nu-dissip}
    \frac{d}{dt}\nu(\phi_t)
=-\int_M (s_{\phi_t}-\bar s)^2\,\omega_{\phi_t}^n
=-g(\phi_t)^2.
\end{align}
 Combining (\ref{absconal}) and 
(\ref{eq:nu-dissip}),
\[
\nu(\alpha(t))-\nu(\alpha(s))
=-\int_s^t g(\alpha(r))^2\,dr
\le -\frac12\int_s^t |\dot\alpha|^2(r)\,dr-\frac12\int_s^t g(\alpha(r))^2\,dr
\]
for all $s,t\in(0,+\infty)$. This yields item $(2)$ of Definition \ref{dgf}, and thereby completes the proof. \hfill\(\square\)

\subsubsection{Proof of Theorem \ref{thmcal}}
Let $g:\mathcal H\to\mathbb R$ be the function defined in \eqref{sugc}. By
Lemma~\ref{lem:slope-Kenergy} we have that
\[
g(\phi)=|\nabla\nu|(\phi)=|\partial\nu|(\phi)\quad\forall \phi\in\mathcal H,
\]
and that $g$ is a strong upper gradient of $\nu$. Moreover, both $\nu$ and $g$ are lower semicontinuous; see Remark~\ref{rcelsc}. The result will follow from Corollary~\ref{cor_gms},  once we verify that
$\phi_0\in\mathcal M(g)$; we do this below.
\smallbreak\noindent
Consider the curve $\alpha:(0,+\infty)\to\mathcal H$ given by $\alpha(t):=\phi_t$. By
Lemma~\ref{lem:gfc}, $\alpha$ is a gradient flow of $\nu$ with respect to the upper gradient
$g$. Since $(t,z)\longmapsto \phi_t(z)$ is smooth on $[0,+\infty)\times M$, we have
\[
    \lim_{t\longrightarrow0^+}\nu\circ\alpha(t)=\lim_{t\longrightarrow0^+}\nu(\phi_t)=\nu(\phi_0)\quad\text{and}\quad \lim_{t\longrightarrow0^+}g\circ\alpha(t)=\lim_{t\longrightarrow0^+}g(\phi_t)=g(\phi_0).
\]
Using the path-length estimate
\eqref{absconal}, we obtain that
\[
d(\phi_0,\alpha(t))=d(\phi_0,\phi_t)\le \int_0^t\Bigl(\int_M \dot\phi_r^{\,2}\,\omega_{\phi_r}^n\Bigr)^{1/2}\,dr
=\int_0^t g(\phi_r)\,dr\quad\forall t\in(0,+\infty).
\]
From this, we see that $d(\alpha(t),\phi_0)\longrightarrow0$ as $t\longrightarrow 0^+$; whence the result follows. \hfill\(\square\)


\appendix
\section{Metric derivatives and strong upper gradients}
In this appendix we present brief preliminaries on metric derivatives and upper gradients. Our presentation follows \cite[Chapter 1]{AGS_2008} and \cite[Section 2]{HM_2019}. 
\subsection{Metric derivatives}
Let $(X,d)$ be a metric space. A \emph{curve} is a continuous map $v:I\to X$
from an interval $I \subseteq \mathbb{R}$ into $X$. In what follows, we restrict attention to curves defined on the interval $(0,+\infty)$. Observe that continuity is purely a topological property. For many purposes, however, more quantitative notions of continuity are needed. 
\begin{definition}
We say that a curve $v:(0,+\infty)\to X$ is locally absolutely continuous if there exists a nonnegative function $m\in L^1_{\operatorname{loc}}(0,+\infty)$ such that 
\begin{align*}
    d\big(v(t), v(s)\big) \le \int_{s}^t m(r)\, dr
\end{align*}
for all $s,t\in(0,+\infty)$ satisfying $s\le t$.
\end{definition}\noindent
The above definition extends the classical notion of locally absolutely continuous functions from $(0,+\infty)$ into $\mathbb{R}^d$. We recall that \(L^1_{\mathrm{loc}}(0,+\infty)\) is the space of measurable functions from 
\((0,+\infty)\) into $\mathbb R$ whose restriction to every compact subset 
\(K \subset (0,+\infty)\) is integrable.

 \begin{definition}
 The metric derivative of a curve $v:(0,+\infty)\to X$ at a point  $t\in(0,+\infty)$ is given by
\[
|\dot v|(t) := \lim_{h \to 0} \frac{d\big(v(t+h),\,v(t)\big)}{|h|},
\]
provided this limit exists.
 \end{definition}
 \noindent
It is a standard fact that, for a locally absolutely continuous curve $v:(0,+\infty)\to X$, its metric derivative exists almost everywhere and $|\dot v|$ belongs to $L^1_{\operatorname{loc}}(0,+\infty)$. 
\begin{remark} If $X$ is a normed space and $v:(0,+\infty)\to X$  a differentiable curve, then 
\[
|\dot v|(t) = \|v'(t)\|_{X}\quad \forall t\in(0,+\infty). 
\]
\end{remark}
 
\subsection{Strong upper gradients}
Let $(X,d)$ be a metric space and $\Phi:X\to \mathbb R\cup\{+\infty\}$ a proper functional. 
\begin{definition}\label{sug}
    We say that a nonnegative function $g:X\to\mathbb R\cup\{+\infty\}$ is a \textit{(strong) upper gradient} of $\Phi$ if for every locally absolutely continuous curve $v:(0,+\infty)\to X$ with $\{v(t):\, t\in(0,+\infty)\}\subseteq\dom \Phi$, the function $g\circ v:(0,+\infty)\to \mathbb R\cup\{+\infty\}$ is Borel measurable and
\begin{align}\label{sugin}
    |\Phi(v(t)) - \Phi(v(s))| \;\le\; \int_s^t g(v(r))\,|\dot v|(r)\,dr
\end{align}
for all $s,t\in(0,+\infty)$ satisfying $s\le t$.
\end{definition}\noindent
Definition \ref{sug} essentially\footnote{In principle, \cite[Definition~1.2.1]{AGS_2008} is ill defined if there exist times $s,t\in(0,+\infty)$ such that $v(s),v(t)\notin\dom\Phi$, since the corresponding expression involves $+\infty-\infty$, which is not defined in the usual extended-real arithmetic.} coincides with the standard one in \cite[Definition~1.2.1]{AGS_2008}, although here it is stated in terms of locally absolutely continuous curves rather than globally absolutely continuous ones. Note that in the definition of upper gradient, the requirement that $g\circ v$ be Borel holds automatically if $g$ is Borel, and that the quantity on the right-hand side of (\ref{sugin}) is allowed to be infinite.

\begin{remark}
If $X$ is a normed space and $\Phi:X\to\mathbb{R}$ a continuously differentiable functional, 
then a Borel function $g:X\to\mathbb{R}\cup\{+\infty\}$ is an upper gradient of $\Phi$ if 
\[
    \|D\Phi(u)\|_{X^*}\le g(u) \quad\forall u\in X.
\]
\end{remark}

\bibliography{references}{}
\bibliographystyle{plain}
\nocite{*}

\end{document}